\def\Date{2008/04/07}


\ifx\pdfoutput\jamaisdefined\else
\input supp-pdf.tex \pdfoutput=1 \pdfcompresslevel=9

\fi

%

\magnification=1200
\hsize=11.25cm
\vsize=18cm
\parskip 0pt
\parindent=12pt
\voffset=1cm
\hoffset=1cm



\catcode'32=9

\font\tenpc=cmcsc10
\font\eightpc=cmcsc8
\font\eightrm=cmr8
\font\eighti=cmmi8
\font\eightsy=cmsy8
\font\eightbf=cmbx8
\font\eighttt=cmtt8
\font\eightit=cmti8
\font\eightsl=cmsl8
\font\sixrm=cmr6
\font\sixi=cmmi6
\font\sixsy=cmsy6
\font\sixbf=cmbx6

\skewchar\eighti='177 \skewchar\sixi='177
\skewchar\eightsy='60 \skewchar\sixsy='60

\catcode`@=11

\def\tenpoint{%
  \textfont0=\tenrm \scriptfont0=\sevenrm \scriptscriptfont0=\fiverm
  \def\rm{\fam\z@\tenrm}%
  \textfont1=\teni \scriptfont1=\seveni \scriptscriptfont1=\fivei
  \def\oldstyle{\fam\@ne\teni}%
  \textfont2=\tensy \scriptfont2=\sevensy \scriptscriptfont2=\fivesy
  \textfont\itfam=\tenit
  \def\it{\fam\itfam\tenit}%
  \textfont\slfam=\tensl
  \def\sl{\fam\slfam\tensl}%
  \textfont\bffam=\tenbf \scriptfont\bffam=\sevenbf
  \scriptscriptfont\bffam=\fivebf
  \def\bf{\fam\bffam\tenbf}%
  \textfont\ttfam=\tentt
  \def\tt{\fam\ttfam\tentt}%
  \abovedisplayskip=12pt plus 3pt minus 9pt
  \abovedisplayshortskip=0pt plus 3pt
  \belowdisplayskip=12pt plus 3pt minus 9pt
  \belowdisplayshortskip=7pt plus 3pt minus 4pt
  \smallskipamount=3pt plus 1pt minus 1pt
  \medskipamount=6pt plus 2pt minus 2pt
  \bigskipamount=12pt plus 4pt minus 4pt
  \normalbaselineskip=12pt
  \setbox\strutbox=\hbox{\vrule height8.5pt depth3.5pt width0pt}%
  \let\bigf@ntpc=\tenrm \let\smallf@ntpc=\sevenrm
  \let\petcap=\tenpc
  \normalbaselines\rm}

\def\eightpoint{%
  \textfont0=\eightrm \scriptfont0=\sixrm \scriptscriptfont0=\fiverm
  \def\rm{\fam\z@\eightrm}%
  \textfont1=\eighti \scriptfont1=\sixi \scriptscriptfont1=\fivei
  \def\oldstyle{\fam\@ne\eighti}%
  \textfont2=\eightsy \scriptfont2=\sixsy \scriptscriptfont2=\fivesy
  \textfont\itfam=\eightit
  \def\it{\fam\itfam\eightit}%
  \textfont\slfam=\eightsl
  \def\sl{\fam\slfam\eightsl}%
  \textfont\bffam=\eightbf \scriptfont\bffam=\sixbf
  \scriptscriptfont\bffam=\fivebf
  \def\bf{\fam\bffam\eightbf}%
  \textfont\ttfam=\eighttt
  \def\tt{\fam\ttfam\eighttt}%
  \abovedisplayskip=9pt plus 2pt minus 6pt
  \abovedisplayshortskip=0pt plus 2pt
  \belowdisplayskip=9pt plus 2pt minus 6pt
  \belowdisplayshortskip=5pt plus 2pt minus 3pt
  \smallskipamount=2pt plus 1pt minus 1pt
  \medskipamount=4pt plus 2pt minus 1pt
  \bigskipamount=9pt plus 3pt minus 3pt
  \normalbaselineskip=9pt
  \setbox\strutbox=\hbox{\vrule height7pt depth2pt width0pt}%
  \let\bigf@ntpc=\eightrm \let\smallf@ntpc=\sixrm
  \let\petcap=\eightpc
  \normalbaselines\rm}
\catcode`@=12

\tenpoint


\long\def\irmaaddress{{%
\bigskip
\eightpoint
\rightline{\quad
\vtop{\halign{\hfil##\hfil\cr
I.R.M.A. UMR 7501\cr
Universit\'e Louis Pasteur et CNRS,\cr
7, rue Ren\'e-Descartes\cr
F-67084 Strasbourg, France\cr
{\tt guoniu@math.u-strasbg.fr}\cr}}\quad}
}}



\catcode`\@=11
\def\pc#1#2|{{\bigf@ntpc #1\penalty \@MM\hskip\z@skip\smallf@ntpc%
	\uppercase{#2}}}
\catcode`\@=12

\def\pointir{\discretionary{.}{}{.\kern.35em---\kern.7em}\nobreak
   \hskip 0em plus .3em minus .4em }

\def\qed{\quad\raise -2pt\hbox{\vrule\vbox to 10pt{\hrule width 4pt
   \vfill\hrule}\vrule}}

\def\rem#1|{\par\medskip{{\it #1}\pointir}}

\def\vspace[#1]{\noalign{\vskip#1}}

\def\abstract#1{\vbox{\eightpoint\narrower\narrower 
\pc ABSTRACT|\pointir #1}}


\def\section#1{\goodbreak\par\vskip .3cm\centerline{\bf #1}
   \par\nobreak\vskip 3pt }

\long\def\th#1|#2\endth{\par\medbreak
   {\petcap #1\pointir}{\it #2}\par\medbreak}

\def\article#1|#2|#3|#4|#5|#6|#7|
    {{\leftskip=7mm\noindent
     \hangindent=2mm\hangafter=1
     \llap{{\tt [#1]}\hskip.35em}{\petcap#2}\pointir
     #3, {\sl #4}, {\bf #5} ({\oldstyle #6}),
     pp.\nobreak\ #7.\par}}
\def\livre#1|#2|#3|#4|
    {{\leftskip=7mm\noindent
    \hangindent=2mm\hangafter=1
    \llap{{\tt [#1]}\hskip.35em}{\petcap#2}\pointir
    {\sl #3}, #4.\par}}
\def\divers#1|#2|#3|
    {{\leftskip=7mm\noindent
    \hangindent=2mm\hangafter=1
     \llap{{\tt [#1]}\hskip.35em}{\petcap#2}\pointir
     #3.\par}}



\catcode`\@=11
\def\c@rr@#1{\vbox{%
  \hrule height \ep@isseur%
   \hbox{\vrule width\ep@isseur\vbox to \t@ille{%
           \vfil\hbox  to \t@ille{\hfil#1\hfil}\vfil}%
            \vrule width\ep@isseur}%
      \hrule height \ep@isseur}}
\def\ytableau#1#2#3#4{\vbox{%
  \gdef\ep@isseur{#2}
   \gdef\t@ille{#1}
    \def\\##1{\c@rr@{$#3 ##1$}}
  \lineskiplimit=-30cm \baselineskip=\t@ille%
    \advance \baselineskip by \ep@isseur%
     \halign{%
      \hfil$##$\hfil&&\kern -\ep@isseur%
       \hfil$##$\hfil \crcr#4\crcr}}}%
\catcode`\@=12

\def\Grille{\setbox13=\vbox to 5mm{\hrule width 110mm\vfill}
\setbox13=\vbox{\offinterlineskip
   \copy13\copy13\copy13\copy13\copy13\copy13\copy13\copy13
   \copy13\copy13\copy13\copy13\box13\hrule width 110mm}
\setbox14=\hbox to 5mm{\vrule height 65mm\hfill}
\setbox14=\hbox{\copy14\copy14\copy14\copy14\copy14\copy14
   \copy14\copy14\copy14\copy14\copy14\copy14\copy14\copy14
   \copy14\copy14\copy14\copy14\copy14\copy14\copy14\copy14\box14}
\ht14=0pt\dp14=0pt\wd14=0pt
\setbox13=\vbox to 0pt{\vss\box13\offinterlineskip\box14}
\wd13=0pt\box13}


\def\fleche(#1,#2)\dir(#3,#4)\long#5{%
\noalign{\nointerlineskip\leftput(#1,#2){\vector(#3,#4){#5}}\nointerlineskip}}


\def\hfl#1#2#3{\smash{\mathop{\hbox to#3{\rightarrowfill}}\limits
^{\scriptstyle#1}_{\scriptstyle#2}}}

\def\gfl#1#2#3{\smash{\mathop{\hbox to#3{\leftarrowfill}}\limits
^{\scriptstyle#1}_{\scriptstyle#2}}}


 \message{`lline' & `vector' macros from LaTeX}
 \catcode`@=11
\def\{{\relax\ifmmode\lbrace\else$\lbrace$\fi}
\def\}{\relax\ifmmode\rbrace\else$\rbrace$\fi}
\def\newcount{\alloc@0\count\countdef\insc@unt}
\def\newdimen{\alloc@1\dimen\dimendef\insc@unt}
\def\newwrite{\alloc@7\write\chardef\sixt@@n}

\newwrite\@unused
\newcount\@tempcnta
\newcount\@tempcntb
\newdimen\@tempdima
\newdimen\@tempdimb
\newbox\@tempboxa

\def\@spaces{\space\space\space\space}
\def\@whilenoop#1{}
\def\@whiledim#1\do #2{\ifdim #1\relax#2\@iwhiledim{#1\relax#2}\fi}
\def\@iwhiledim#1{\ifdim #1\let\@nextwhile=\@iwhiledim
        \else\let\@nextwhile=\@whilenoop\fi\@nextwhile{#1}}
\def\@badlinearg{\@latexerr{Bad \string\line\space or \string\vector
   \space argument}}
\def\@latexerr#1#2{\begingroup
\edef\@tempc{#2}\expandafter\errhelp\expandafter{\@tempc}%
\def\@eha{Your command was ignored.
^^JType \space I <command> <return> \space to replace it
  with another command,^^Jor \space <return> \space to continue without it.}
\def\@ehb{You've lost some text. \space \@ehc}
\def\@ehc{Try typing \space <return>
  \space to proceed.^^JIf that doesn't work, type \space X <return> \space to
  quit.}
\def\@ehd{You're in trouble here.  \space\@ehc}

\typeout{LaTeX error. \space See LaTeX manual for explanation.^^J
 \space\@spaces\@spaces\@spaces Type \space H <return> \space for
 immediate help.}\errmessage{#1}\endgroup}
\def\typeout#1{{\let\protect\string\immediate\write\@unused{#1}}}

\font\tenln    = line10
\font\tenlnw   = linew10

\newdimen\@wholewidth
\newdimen\@halfwidth
\newdimen\unitlength 

\unitlength =1pt


\def\thinlines{\let\@linefnt\tenln \let\@circlefnt\tencirc
  \@wholewidth\fontdimen8\tenln \@halfwidth .5\@wholewidth}
\def\thicklines{\let\@linefnt\tenlnw \let\@circlefnt\tencircw
  \@wholewidth\fontdimen8\tenlnw \@halfwidth .5\@wholewidth}

\def\linethickness#1{\@wholewidth #1\relax \@halfwidth .5\@wholewidth}

\newif\if@negarg

\def\lline(#1,#2)#3{\@xarg #1\relax \@yarg #2\relax
\@linelen=#3\unitlength
\ifnum\@xarg =0 \@vline
  \else \ifnum\@yarg =0 \@hline \else \@sline\fi
\fi}

\def\@sline{\ifnum\@xarg< 0 \@negargtrue \@xarg -\@xarg \@yyarg -\@yarg
  \else \@negargfalse \@yyarg \@yarg \fi
\ifnum \@yyarg >0 \@tempcnta\@yyarg \else \@tempcnta -\@yyarg \fi
\ifnum\@tempcnta>6 \@badlinearg\@tempcnta0 \fi
\setbox\@linechar\hbox{\@linefnt\@getlinechar(\@xarg,\@yyarg)}%
\ifnum \@yarg >0 \let\@upordown\raise \@clnht\z@
   \else\let\@upordown\lower \@clnht \ht\@linechar\fi
\@clnwd=\wd\@linechar
\if@negarg \hskip -\wd\@linechar \def\@tempa{\hskip -2\wd\@linechar}\else
     \let\@tempa\relax \fi
\@whiledim \@clnwd <\@linelen \do
  {\@upordown\@clnht\copy\@linechar
   \@tempa
   \advance\@clnht \ht\@linechar
   \advance\@clnwd \wd\@linechar}%
\advance\@clnht -\ht\@linechar
\advance\@clnwd -\wd\@linechar
\@tempdima\@linelen\advance\@tempdima -\@clnwd
\@tempdimb\@tempdima\advance\@tempdimb -\wd\@linechar
\if@negarg \hskip -\@tempdimb \else \hskip \@tempdimb \fi
\multiply\@tempdima \@m
\@tempcnta \@tempdima \@tempdima \wd\@linechar \divide\@tempcnta \@tempdima
\@tempdima \ht\@linechar \multiply\@tempdima \@tempcnta
\divide\@tempdima \@m
\advance\@clnht \@tempdima
\ifdim \@linelen <\wd\@linechar
   \hskip \wd\@linechar
  \else\@upordown\@clnht\copy\@linechar\fi}

\def\@hline{\ifnum \@xarg <0 \hskip -\@linelen \fi
\vrule height \@halfwidth depth \@halfwidth width \@linelen
\ifnum \@xarg <0 \hskip -\@linelen \fi}

\def\@getlinechar(#1,#2){\@tempcnta#1\relax\multiply\@tempcnta 8
\advance\@tempcnta -9 \ifnum #2>0 \advance\@tempcnta #2\relax\else
\advance\@tempcnta -#2\relax\advance\@tempcnta 64 \fi
\char\@tempcnta}

\def\vector(#1,#2)#3{\@xarg #1\relax \@yarg #2\relax
\@linelen=#3\unitlength
\ifnum\@xarg =0 \@vvector
  \else \ifnum\@yarg =0 \@hvector \else \@svector\fi
\fi}

\def\@hvector{\@hline\hbox to 0pt{\@linefnt
\ifnum \@xarg <0 \@getlarrow(1,0)\hss\else
    \hss\@getrarrow(1,0)\fi}}

\def\@vvector{\ifnum \@yarg <0 \@downvector \else \@upvector \fi}

\def\@svector{\@sline
\@tempcnta\@yarg \ifnum\@tempcnta <0 \@tempcnta=-\@tempcnta\fi
\ifnum\@tempcnta <5
  \hskip -\wd\@linechar
  \@upordown\@clnht \hbox{\@linefnt  \if@negarg
  \@getlarrow(\@xarg,\@yyarg) \else \@getrarrow(\@xarg,\@yyarg) \fi}%
\else\@badlinearg\fi}

\def\@getlarrow(#1,#2){\ifnum #2 =\z@ \@tempcnta='33\else
\@tempcnta=#1\relax\multiply\@tempcnta \sixt@@n \advance\@tempcnta
-9 \@tempcntb=#2\relax\multiply\@tempcntb \tw@
\ifnum \@tempcntb >0 \advance\@tempcnta \@tempcntb\relax
\else\advance\@tempcnta -\@tempcntb\advance\@tempcnta 64
\fi\fi\char\@tempcnta}

\def\@getrarrow(#1,#2){\@tempcntb=#2\relax
\ifnum\@tempcntb < 0 \@tempcntb=-\@tempcntb\relax\fi
\ifcase \@tempcntb\relax \@tempcnta='55 \or
\ifnum #1<3 \@tempcnta=#1\relax\multiply\@tempcnta
24 \advance\@tempcnta -6 \else \ifnum #1=3 \@tempcnta=49
\else\@tempcnta=58 \fi\fi\or
\ifnum #1<3 \@tempcnta=#1\relax\multiply\@tempcnta
24 \advance\@tempcnta -3 \else \@tempcnta=51\fi\or
\@tempcnta=#1\relax\multiply\@tempcnta
\sixt@@n \advance\@tempcnta -\tw@ \else
\@tempcnta=#1\relax\multiply\@tempcnta
\sixt@@n \advance\@tempcnta 7 \fi\ifnum #2<0 \advance\@tempcnta 64 \fi
\char\@tempcnta}

\def\@vline{\ifnum \@yarg <0 \@downline \else \@upline\fi}

\def\@upline{\hbox to \z@{\hskip -\@halfwidth \vrule
  width \@wholewidth height \@linelen depth \z@\hss}}

\def\@downline{\hbox to \z@{\hskip -\@halfwidth \vrule
  width \@wholewidth height \z@ depth \@linelen \hss}}

\def\@upvector{\@upline\setbox\@tempboxa\hbox{\@linefnt\char'66}\raise
     \@linelen \hbox to\z@{\lower \ht\@tempboxa\box\@tempboxa\hss}}

\def\@downvector{\@downline\lower \@linelen
      \hbox to \z@{\@linefnt\char'77\hss}}

\thinlines

\newcount\@xarg
\newcount\@yarg
\newcount\@yyarg
\newcount\@multicnt
\newdimen\@xdim
\newdimen\@ydim
\newbox\@linechar
\newdimen\@linelen
\newdimen\@clnwd
\newdimen\@clnht
\newdimen\@dashdim
\newbox\@dashbox
\newcount\@dashcnt
 \catcode`@=12


\newbox\tbox
\newbox\tboxa

\def\leftzer#1{\setbox\tbox=\hbox to 0pt{#1\hss}%
     \ht\tbox=0pt \dp\tbox=0pt \box\tbox}

\def\rightzer#1{\setbox\tbox=\hbox to 0pt{\hss #1}%
     \ht\tbox=0pt \dp\tbox=0pt \box\tbox}

\def\centerzer#1{\setbox\tbox=\hbox to 0pt{\hss #1\hss}%
     \ht\tbox=0pt \dp\tbox=0pt \box\tbox}

%
\def\image(#1,#2)#3{\vbox to #1{\offinterlineskip
    \vss #3 \vskip #2}}


\def\leftput(#1,#2)#3{\setbox\tboxa=\hbox{%
    \kern #1\unitlength
    \raise #2\unitlength\hbox{\leftzer{#3}}}%
    \ht\tboxa=0pt \wd\tboxa=0pt \dp\tboxa=0pt\box\tboxa}

\def\rightput(#1,#2)#3{\setbox\tboxa=\hbox{%
    \kern #1\unitlength
    \raise #2\unitlength\hbox{\rightzer{#3}}}%
    \ht\tboxa=0pt \wd\tboxa=0pt \dp\tboxa=0pt\box\tboxa}

\def\centerput(#1,#2)#3{\setbox\tboxa=\hbox{%
    \kern #1\unitlength
    \raise #2\unitlength\hbox{\centerzer{#3}}}%
    \ht\tboxa=0pt \wd\tboxa=0pt \dp\tboxa=0pt\box\tboxa}

\unitlength=1mm

\def\cput(#1,#2)#3{\noalign{\nointerlineskip\centerput(#1,#2){#3}
                             \nointerlineskip}}


\ifx\pdfoutput\jamaisdefined
\input epsf

\fi


\parskip 0pt plus 1pt

\def\article#1|#2|#3|#4|#5|#6|#7|
    {{\leftskip=7mm\noindent
     \hangindent=2mm\hangafter=1
     \llap{{\tt [#1]}\hskip.35em}{#2},\quad %
     #3, {\sl #4}, {\bf #5} ({\oldstyle #6}),
     pp.\nobreak\ #7.\par}}
\def\livre#1|#2|#3|#4|
    {{\leftskip=7mm\noindent
    \hangindent=2mm\hangafter=1
    \llap{{\tt [#1]}\hskip.35em}{#2},\quad %
    {\sl #3}, #4.\par}}
\def\divers#1|#2|#3|
    {{\leftskip=7mm\noindent
    \hangindent=2mm\hangafter=1
     \llap{{\tt [#1]}\hskip.35em}{#2},\quad %
     #3.\par}}


%

\def\setB{\mathop{\cal B}}
\def\setT{\mathop{\cal T}}
\def\setF{\mathop{\cal F}}


\rightline{\Date}
\bigskip

\centerline{\bf Yet another generalization of }
\centerline{\bf Postnikov's hook length formula for binary trees}
\bigskip
\centerline{Guo-Niu HAN}
\bigskip\medskip

\abstract{
We discover another one-parameter generalization of Postnikov's hook 
length formula for binary trees. 
The particularity of our formula is that the hook length $h_v$ appears as 
an exponent.  As an application, we derive another simple hook length formula 
for binary trees when the underlying parameter takes the value~$1/2$.
}

\bigskip

\def\sec{1}
\section{\sec. Introduction} 
Consider the set $\setB(n)$ of all binary trees with $n$ vertices. 
For each vertex~$v$ of $T\in \setB(n)$ the {\it hook length} 
of~$v$, denoted by $h_v$, or just $h$ for short, is 
the number of descendants of $v$ (including $v$).
The following hook length formula for binary trees
$$
\sum_{T\in\setB(n)}  \prod_{v\in T} \bigl(1+{1\over h_v}\bigr)
= {2^n\over n!}(n+1)^{n-1}
\leqno{(1)}
$$
was discovered by Postnikov [Po04].  Further combinatorial proofs and 
extensions have been proposed by several authors [CY08, GS06, MY07, Se08]. 
In particular, Lascoux conjectured the following 
one-parameter generalization 
$$
\sum_{T\in\setB(n)} \prod_{v\in T} \bigl(x+{1\over h_v}\bigr)
= {1\over (n+1)!}  \prod_{k=0}^{n-1} ((n+1+k)x+n+1-k),
\leqno{(2)}
$$
which was subsequently proved by Du-Liu [DL08].
The latter generalization appears to be very natural, because
the {\it left-hand side} of (2) can be obtained from the left-hand side of (1)
by replacing~1 by $x$.

\medskip
It is also natural to look for an extension of (1) by introducing a new 
variable $z$ in the {\it right-hand side}, namely by replacing 
$2^n(n+1)^{n-1}/n!$ by $2^nz(n+z)^{n-1}/n!$  It so happens that the 
corresponding left-hand side is also a sum on binary trees, but this time the
hook length $h_v$ appears as an exponent. The purpose of this Note is 
to prove the following Theorem.

\proclaim Theorem 1.
For each positive integer $n$ we have
$$
\sum_{T\in\setB(n)}  \prod_{v\in T} {(z+h)^{h-1}\over h (2z+h-1)^{h-2}}
= {2^nz\over n!}(n+z)^{n-1}.
\leqno{(3)}
$$

With $z=1$ in (3) we recover Postnikov's 
identity (1). The following corollary is derived from our identity (3)
by taking $z=1/2$.

\proclaim Corollary 2.
For each positive integer $n$ we have
$$
\sum_{T\in\setB(n)}  \prod_{v\in T} \bigl(1+{1\over  2h}\bigr)^{h-1}
= {  (2n+1)^{n-1} \over n!}.
\leqno{(4)}
$$


\def\sec{2}
\section{\sec. Proof of the Theorem} 

Let us take an example before proving the Theorem.  There are five 
binary trees with $n=3$ vertices:

%

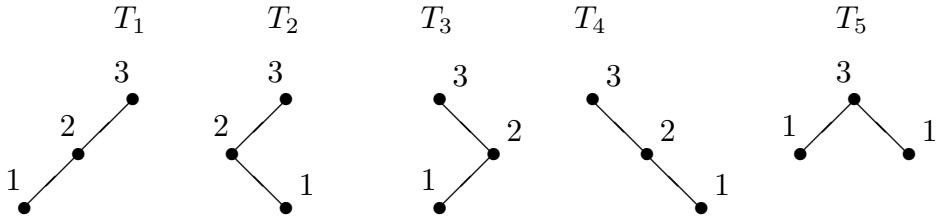
\long\def\maplebegin#1\mapleend{}

\maplebegin

# --------------- begin maple ----------------------

# Copy the following text  to "makefig.mpl"
# then in maple > read("makefig.mpl");
# it will create a file "z_fig_by_maple.tex"

#\unitlength=1pt

Hu:= 6; # height quantities

X0:=5.0; Y0:=0; # origin position

File:=fopen("z_fig_by_maple.tex", WRITE);

pline:=proc(x,y) # X0,Y0 = offset
local a,b,len;
	len:=1;
	fprintf(File, "\\pline(
end;

nline:=proc(x,y) # X0,Y0 = offset
local a,b,len;
	len:=1;
	fprintf(File, "\\nline(
end;

mydot:=proc(x,y) # X0,Y0 = offset
local a,b,len;
	len:=1;
	fprintf(File, "\\mydot(
end;

mylabel:=proc(x,y, text) # X0,Y0 = offset
local a,b,len;
	len:=1;
	fprintf(File, "\\mylabel(
end;

mylabel2:=proc(x,y, text) # X0,Y0 = offset
local a,b,len;
	len:=1;
	fprintf(File, "\\mylabel(
end;

dotlabel:=proc(x,y,t) mydot(x,y); mylabel(x,y,t); end;
dotlabel2:=proc(x,y,t) mydot(x,y); mylabel2(x,y,t); end;

DXX:=23;

X0:=-2.6*DXX;
pline(2,1); pline(3,2); 
dotlabel(2,1, "1"); dotlabel(3,2, "2"); dotlabel(4,3, "3");
mylabel(4,4, "$T_1$");

X0:=X0+DXX;

nline(2,2); pline(2,2); 
dotlabel2(3,1, "1"); dotlabel(2,2, "2"); dotlabel(3,3, "3");
mylabel(3,4, "$T_2$");

X0:=X0+DXX;
nline(2,3); pline(2,1); 
dotlabel(2,1, "1"); dotlabel2(3,2, "2"); dotlabel2(2,3, "3");
mylabel(2,4, "$T_3$");

X0:=X0+DXX; 
nline(1,3); nline(2,2); 
dotlabel2(3,1, "1"); dotlabel2(2,2, "2"); dotlabel2(1,3, "3");
mylabel(1,4, "$T_4$");

X0:=X0+DXX;
pline(1,2); nline(2,3); 
dotlabel(1,2, "1"); dotlabel2(3,2, "1"); dotlabel(2,3, "3");
mylabel(2,4, "$T_5$");

fclose(File);

# -------------------- end maple -------------------------

\mapleend


\newbox\boxarbre
\def\pline(#1,#2)#3|{\leftput(#1,#2){\lline(1,1){#3}}}
\def\nline(#1,#2)#3|{\leftput(#1,#2){\lline(1,-1){#3}}}
\def\mydot(#1,#2)|{\leftput(#1,#2){$\bullet$}}
\def\mylabel(#1,#2)#3|{\leftput(#1,#2){#3}}
\setbox\boxarbre=\vbox{\vskip
30mm\offinterlineskip 
%
\pline(-47.8,6.0)6.0|
\pline(-41.8,12.0)6.0|
\mydot(-48.6,5.2)|
\mylabel(-49.8,8.0){1}|
\mydot(-42.6,11.2)|
\mylabel(-43.8,14.0){2}|
\mydot(-36.6,17.2)|
\mylabel(-37.8,20.0){3}|
\mylabel(-37.8,26.0){$T_1$}|
\nline(-24.8,12.0)6.0|
\pline(-24.8,12.0)6.0|
\mydot(-19.6,5.2)|
\mylabel(-17.3,7.5){1}|
\mydot(-25.6,11.2)|
\mylabel(-26.8,14.0){2}|
\mydot(-19.6,17.2)|
\mylabel(-20.8,20.0){3}|
\mylabel(-20.8,26.0){$T_2$}|
\nline(-1.8,18.0)6.0|
\pline(-1.8,6.0)6.0|
\mydot(-2.6,5.2)|
\mylabel(-3.8,8.0){1}|
\mydot(3.4,11.2)|
\mylabel(5.7,13.5){2}|
\mydot(-2.6,17.2)|
\mylabel(-0.3,19.5){3}|
\mylabel(-3.8,26.0){$T_3$}|
\nline(15.2,18.0)6.0|
\nline(21.2,12.0)6.0|
\mydot(26.4,5.2)|
\mylabel(28.7,7.5){1}|
\mydot(20.4,11.2)|
\mylabel(22.7,13.5){2}|
\mydot(14.4,17.2)|
\mylabel(16.7,19.5){3}|
\mylabel(13.2,26.0){$T_4$}|
\pline(38.2,12.0)6.0|
\nline(44.2,18.0)6.0|
\mydot(37.4,11.2)|
\mylabel(36.2,14.0){1}|
\mydot(49.4,11.2)|
\mylabel(51.7,13.5){1}|
\mydot(43.4,17.2)|
\mylabel(42.2,20.0){3}|
\mylabel(42.2,26.0){$T_5$}|%
}
$$
\kern-4mm\box\boxarbre
$$
\noindent
The hook lengths of $T_1, T_2, T_3, T_4$ are all the same $1,2,3$; but
the hook lengths of $T_5$ are $1,1,3$. The left-hand side of (3) 
is then equal to
$$
4\times{1 \over (2z)^{-1}}\cdot
{ (z+2)^1 \over 2} \cdot
{ (z+3)^2 \over 3(2z+1) }
+
{1 \over (2z)^{-1}}\cdot
{1 \over (2z)^{-1}}\cdot
{ (z+3)^2 \over 3(2z+1) }
={2^3z(z+3)^2\over 3!}.
$$
\medskip
Let $y(x)$ be a formal power series in $x$ such that
$$
y(x)=e^{xy(x)}. \leqno{(5)}
$$
By the Lagrange inversion formula $y(x)^z$ has the following explicit 
expansion:
$$
y(x)^{z}=\sum_{n\geq 0} z(n+z)^{n-1} {x^n\over n!}. \leqno{(6)}
$$
Since $y^{2z}= (y^z)^2$ we have
$$
\sum_{n\geq 0} 2z(n+2z)^{n-1} {x^n\over n!}
=\bigl(\sum_{n\geq 0} z(n+z)^{n-1} {x^n\over n!}\bigr)^2. \leqno{(7)}
$$
Comparing the coefficients of $x^n$ on both sides of (7) yields the 
following Lemma.

\proclaim Lemma 3.
We have
$$
{2z (n+2z)^{n-1}\over n!}
=
\sum_{k=0}^n
{z(k+z)^{k-1} \over k!} \times
{z(n-k+z)^{n-k-1} \over (n-k)!}.\leqno{(8)}
$$

{\it Proof of the Theorem}.
Let 
$$
P(n)= 
\sum_{T\in\setB(n)}  \prod_{v\in T} {(z+h)^{h-1}\over h (2z+h-1)^{h-2}}.
$$
With each binary tree $T\in\setB(n)$ ($n\geq 1$) we can associate a 
triplet $(T', T'', u)$, where $T'\in\setB(k)$ ($0\leq k\leq n-1$),
$T''\in\setB(n-1-k)$  and $u$ is a vertex of hook length $h_u=n$.  Hence
$$
P(n)=\sum_{k=0}^{n-1} P(k) P(n-1-k)\times 
{(z+n)^{n-1}\over n (2z+n-1)^{n-2}}.
\leqno{(9)}
$$
It is routine to verify that $P(n)=2^nz(z+n)^{n-1}/n!$ for $n=1,2,3$. 
Suppose that
$P(k)=2^kz(z+k)^{k-1}/k!$ for $k\leq n-1$. 
From identity (9) and Lemma~3 we have 
$$
\leqalignno{
P(n)
&=\sum_{k=0}^{n-1} 
{2^kz(z+k)^{k-1}\over k!}\times
{2^{n-k-1}z(z+n-k-1)^{n-k-2}\over (n-k-1)!}\cr
&\kern 5cm \times {(z+n)^{n-1}\over n (2z+n-1)^{n-2}}\cr
&= {2^nz \over n!}(z+n)^{n-1}.\cr
}
$$
By induction, formula (3) is true for any positive integer $n$. 
\qed


\def\sec{3}
\section{\sec. Concluding and Remarks} 
The right-hand sides of (3) and (4) have been studied by
other authors [GS06, DL08, MY07], but our formula has 
the following two major differences:  (i) the hook length $h_v$
appears as an exponent; (ii) the underlying set remains the 
set of binary trees,  whereas in the above mentioned papers
the summation has been changed to the 
set of $m$-ary trees or plane forests. It is interesting to
compare Corollary~2 with the following results obtained by Du and Liu [DL08].
Note that the right-hand sides of formulas (4), (10) and (11) are all 
identical!

\proclaim Proposition 4.
For each positive integer $n$ we have
$$
\sum_{T\in\setT(n)}  \prod_{v\in I(T)} \bigl({2\over 3}+{1\over  3h}\bigr)
= {  (2n+1)^{n-1} \over n!}, \leqno{(10)}
$$
where $\setT(n)$ is the set of all $3$-ary trees with $n$ internal vertices
and $I(T)$ is the set of all internal vertices of $T$.

\proclaim Proposition 5.
For each positive integer $n$ we have
$$
\sum_{T\in\setF(n)}  \prod_{v\in (T)} \bigl(2-{1\over  h}\bigr)
= {  (2n+1)^{n-1} \over n!}, \leqno{(11)}
$$
where $\setF(n)$ is the set of all plane forests with $n$ vertices.


\vskip 10mm 
\medskip

\bigskip \bigskip


\centerline{References}

{\eightpoint

\bigskip 
\bigskip

\divers CY08|Chen, William Y.C.; Yang, Laura L.M.|On Postnikov's hook length
formula for binary trees, {\sl European Journal of Combinatorics}, 
in press, {\oldstyle 2008}|

\article DL08|Du, Rosena R. X.; Liu, Fu|%
$(k,m)$-Catalan Numbers and Hook Length Polynomials for Plane Trees|
European J. Combin|28|2007|1312--1321|

\divers GS06|Gessel, Ira M.; Seo, Seunghyun|%
A refinement of Cayley's formula for trees,
{\it arXiv:math.CO/0507497}|


\article MY07|Moon, J. W.; Yang, Laura L. M.|%
Postnikov identities and Seo's formulas|Bull. Inst. Combin. Appl.|%
49|2007|21--31|

\divers Po04|Postnikov, Alexander|\quad Permutohedra, associahedra, and beyond,
{\it arXiv:math. CO/0507163}, {\oldstyle 2004}|

\divers Se08|Seo, Seunghyun|%
A combinatorial proof of Postnikov's identity and a generalized enumeration 
of labeled trees,
{\it arXiv:math.CO/0409323}|



\bigskip

\irmaaddress
}
\vfill\eject

\end